\documentclass{article}

\usepackage{amsmath, amssymb, amsfonts, amsthm, url}

\newcommand{\comment}[1]{}

\newenvironment{equations}{\setlength\arraycolsep{1pt}\begin{eqnarray*}}{\end{eqnarray*}}

\newcounter{temp}

\theoremstyle{definition}

\newcommand{\metric}[2]{\langle #1,#2 \rangle}

\newcommand{\ZZ}{\mathbb{Z}}

\title{The Ricci flow does not preserve \\ the set of Zoll metrics}
\author{Dan Jane}
\date{\today}

\begin{document} 

\maketitle

\begin{abstract}  The question of whether or not the set of Zoll metrics on the 2-sphere is connected is still open. Here we show that a naive application of the Ricci flow is not sufficient to answer this problem. \end{abstract}

\section*{Introduction}
The geodesics on the round 2-sphere are all simple, closed curves of the same length. A natural question is to find other surfaces with this property, and was first considered by Jean-Gaston Darboux in the 1880's. In 1903 Otto Zoll demonstrated analytic examples which were surfaces of revolution \cite{zoll-1903-}, and his name became attatched to metrics whose geodesics are all simple, closed curves of the same length.

The \emph{Zoll} surfaces of revolution were shown to be in one-to-one correspondence with odd functions of $(-1, 1) \to (-1, 1)$ by Ren\'e Michel in 1978 \cite{michel-1973-}. More general Zoll metrics were shown to exist by Victor Guillemin; any odd function on the round 2-sphere can be realised as the derivative of a one parameter family of Zoll metrics starting at the round metric.

Basic references for this topic, as well as a more detailed account of the history, can be found in \cite{besse-1978-}; see also \cite[pp. 436 and 488]{berger-2004-}.

An outstanding open problem is the following \cite[Question 200]{berger-2004-}

\medskip

\noindent{\bf Question.} {\it Is the set of Zoll metrics on the 2-sphere connected?}

\medskip

That is, given a Zoll surface $(S^2, g)$ can one give a continuous family of Zoll metrics $\{g_t\}$, with $t \in [0, T]$, such that $g_0 = g$ and $g_T$ is the round metric?

A deformation with the correct endpoints - and a particularly fashionable deformation at this time - is the \emph{Ricci flow}. Introduced by Richard Hamilton in 1982 \cite{hamilton-1982-} it can be loosely described as a non-linear heat flow of the curvature of a metric. In the special case of a surface, up to a space-time reparameterisation in order that the volume is constant, the evolution of the metric is governed by
\begin{displaymath}
\frac{\partial g_t}{\partial t} = -2 \big( K - \overline{K} \big) g_t.
\end{displaymath}
The Gaussian curvature is denoted $K$, a function of $t$, and the bar indicates an average value over the manifold. For further background reading see \cite{ma-2004-}. Notice that this evolution is properly called the \emph{normalised Ricci flow}, but in this paper the distinction is unnecessary.

Bennett Chow showed that, whatever the initial metric on $S^2$, after a finite time along the Ricci flow the curvature is everywhere positive \cite{chow-1991-}. In turn, a result of Hamilton states that if the curvature on a surface is everywhere positive then we have convergence to the round sphere \cite{hamilton-1988-}. Hence given a Zoll metric on the 2-sphere, the Ricci flow gives a natural deformation to the round metric.

The Ricci flow preserves metric symmetries \cite{chow-2004-} - the group of isometries and the set of killing fields. For Zoll surfaces, given the global symmetries of the periodic structure in the tangent bundle, one might hope that the Zoll structure was preserved under the Ricci flow.

\medskip

\noindent{\bf Theorem.} {\it The set of Zoll metrics is not preserved by the Ricci flow.}

\medskip

This result in the tangent bundle echoes the way in which the Ricci flow breaks symplectic symmetries \cite{jane-2007-}.

\medskip

Thanks are due to Martin Kerin for some helpful comments on the overall presentation of the note.

\section*{Proof of theorem}
The argument here is to consider the first approximation of the Ricci flow at some initial metric. More precisely, the proof proceeds by identifying some attribute of Zoll manifolds and proving that, even to the first approximation, along the Ricci flow this quantity no longer reflects a Zoll metric.

A Zoll manifold with common period $2 \pi L$ is called a $C_L$-manifold. The following result of Alan Weinstein is crucial,

\medskip

\noindent{\bf Theorem \cite{weinstein-1974-}.} {\it If $(M, g)$ is an $n$-dimensional $C_L$-manifold, then the ratio
\begin{displaymath}
i(M, g) = \frac{\mathrm{vol}(M, g)}{L^n \,\, \mathrm{vol}(S^n, \mathrm{can})}
\end{displaymath}
is an integer, now called \emph{Weinstein's integer}.}

\medskip

We are considering a volume normalised continuous family $\{g_t\}$, which we assume are all Zoll metrics. We further restrict to surfaces of revolution, so that the map 
\begin{equation} \label{CLtoL}
g_t \mapsto \mbox{ common period of } g_t
\end{equation}
is continuous (consider a geodesic through the North pole). The result of Weinstein implies that this map is onto a discrete set. Thus the common period of $g_0$ is the same as the common period of $g_t$ for all $t$. Consider the length $l(t)$ of a closed curve $\gamma_t$ which is a geodesic for the metric $g_t$. For a particular Zoll metric we will show that the derivative $l'(0)$ along the Ricci flow is non-zero, to give a contradiction.

In fact the very first part of Weinstein's proof shows that the range of the map (\ref{CLtoL}) is discrete. One cannot resist repeating this argument due to its simplicity and beauty, especially here in the case of a $C_L$-surface of volume $4 \pi$;

\medskip

{\it Weinstein's argument. } The geodesic flow of a Zoll metric on $M$ generates a free $S^1$ action of the unit sphere bundle $SM$. The geodesic flow is generated by a vector field $G$, and the dual 1-form is the canonical 1-form $\alpha \in \Omega^1(TM)$ \cite{paternain-1999-}. Thus $d \alpha$ is the curvature 2-form for the induced principal $S^1$-bundle, say $p: SM \to CM = SM/ S^1$.

As the curvature 2-form of a principal bundle, there is a unique 2-form $\Omega$ on $CM$ such that $p^*\Omega = d \alpha / (2 \pi L)$. Normalised in this way $\Omega$ can be identified with the Euler class $e(p)$ of the bundle. Thus if $\mathcal{C}$ is the fundamental 2-cycle of $CM$ then
\begin{displaymath}
\int_{CM} \Omega = \metric{[e(p)]}{\mathcal{C}}.
\end{displaymath}
The right hand side is an integer, which implies $\int_{CM} \Omega$ is also an integer.

As a circle bundle over $M$, the volume of the sphere bundle $SM$ is \begin{displaymath}
\mathrm{vol}(SM) = 2\pi \, \mathrm{vol}(M, g) = 8 \pi^2.
\end{displaymath} We have outlined another bundle $p: SM \to CM$ and so, using the Fubini theorem for fibrations,
\begin{equations}
\mathrm{vol}(SM) &=& \int_{SM} \alpha \wedge d\alpha \\
&=& 2 \pi L \int_{SM} \alpha \wedge p^*(\Omega) \\
&=& 2 \pi L \int_{x \in CM} \left[ \int_{p^{-1}(x)} \alpha \right] \Omega.
\end{equations}

But $\int_{p^{-1}(x)} \alpha = 2 \pi L$ since the surface is $C_L$, and hence
\begin{displaymath}
L^{-2} = \frac{2 \pi \cdot 2 \pi}{8 \pi^2} \int_{CM} \Omega = \frac{1}{2} \int_{CM} \Omega \in \mbox{$\frac{1}{2}$} \ZZ,
\end{displaymath}
a discrete set. \qed

\medskip

Now we outline the Zoll metric with the properties we require in order that $l'(0)$ is non-zero. Gambier's \emph{gong} $\Sigma$ \cite[p. 95]{gambier-1925-} is a Zoll surface of revolution with meridian curve
\begin{displaymath}
r(z) = 4 \bigg( \sqrt{1 + \sqrt{ 1 - z^2 } } - 1 \bigg).
\end{displaymath}
Besides the rotational symmetry, $\Sigma$ also has a reflectional symmetry in the $z = 0$ plane, which implies that the intersection of $\Sigma$ with the $z = 0$ plane is a geodesic. This is the curve we denote by $\gamma_0$, the geodesic of interest for the initial metric $(S^2, g_0) = \Sigma$. 

Further, this reflectional symmetry will be preserved along the Ricci flow: hence $\gamma_t := \gamma_0$ is a geodesic of $g_t$ for all $t$.

Rescale in order that $\Sigma$ has volume $4 \pi$. This implies that $r(0) = 1$ and the period of every geodesic is $2 \pi$. It also implies that the average curvature $\overline{K} = 1$.

We assume $\gamma_0 : [0, 2 \pi] \to S^2$ has constant speed, in order that it is a unit speed geodesic. Now we differentiate the length of $\gamma_t$ in the direction of the Ricci flow;
\begin{displaymath}
l'(0) = \frac{\partial}{\partial t} \Big|_{t=0} l(t) = \frac{\partial}{\partial t} \Big|_{t=0} \int_0^{2\pi} \sqrt{g_t\big( \dot{\gamma_t}(s), \dot{\gamma_t}(s) \big)} \, \, ds,
\end{displaymath}
but $\gamma_t = \gamma_0$ is not a function of $t$, and so
\begin{displaymath}
l'(0) = \int_0^{2\pi} \frac{ \frac{\partial}{\partial t} \Big|_{t=0} g_t\big( \dot{\gamma_0}(s), \dot{\gamma_0}(s) \big)}{2 \sqrt{g_0 \big( \dot{\gamma_0}(s), \dot{\gamma_0}(s) \big)}} \, \, ds.
\end{displaymath}
Since we chose $\gamma_0$ to be a unit speed geodesic $g_0 \big( \dot{\gamma_0}(s), \dot{\gamma_0}(s) \big) = 1$ in the denominator. We use the Ricci flow evolution equation in the numerator,
\begin{displaymath}
l'(0) = \int_0^{2\pi} - \Big( K(\gamma_0(s)) - \overline{K} \Big) g_0 \big( \dot{\gamma_0}(s), \dot{\gamma_0}(s) \big) \, ds,
\end{displaymath}
and again use the observation that the geodesic is unit speed to arrive at
\begin{displaymath}
l'(0)= - \int_0^{2\pi} \big( K(\gamma_0(s)) - \overline{K} \big) \, ds.
\end{displaymath}
Finally, because of the rotational symmetry the Gaussian curvature is constant along $\gamma_0$, say $\kappa$. All that is left to show is that $\kappa$ is not equal to $\overline{K} = 1$.

A calculation shows that for a surface of revolution with parametrisation
\begin{displaymath}
\sigma(z, \theta) = \Big( r(z) \cos(\theta), \, r(z) \sin(\theta), \, z \Big)
\end{displaymath}
the curvature is given by
\begin{displaymath}
K(z, \theta) = K(z) = \frac{- r''(z)}{r(z) \Big( r'(z)^2 + 1 \Big)^2 }.
\end{displaymath}
Since we rescaled the gong in order to normalise the volume, the meridian is now
\begin{displaymath}
r(z) = \Big( \sqrt{1 + \sqrt{ 1 - 16 (cz)^2 } } - 1 \Big) / c,
\end{displaymath}
where $c = \sqrt{2}-1$. Subsituting in gives
\begin{displaymath}
\kappa = K(0) = 4 (2 - \sqrt{2}) \neq 1.
\end{displaymath}

\medskip

\noindent{\bf Remark.} A similar analysis of the lengths of the three simple, closed, plane geodesics of an ellipsoid can be used to show that the Ricci flow does not preserve the class of ellipsoids \cite[Section 3.4]{jane-2008t}.

\bibliography{New_Master}
\bibliographystyle{plain}

\end{document}